\documentclass[12pt]{amsart}
\usepackage{amsmath}

\theoremstyle{plain}
\newtheorem{thm}{Theorem}[section]
\newtheorem{theorem}[thm]{Theorem}
\newtheorem{corollary}[thm]{Corollary}
\newtheorem{lemma}[thm]{Lemma}
\newtheorem{proposition}[thm]{Proposition}
\newtheorem{remark}[thm]{Remark}

\numberwithin{equation}{section}
\newcommand\dd{\frac{d^2}{dx^2}}
\newcommand\ov{\overline}
\newcommand\ti{\widetilde}
\newcommand\lan{\langle}
\newcommand\ran{\rangle}

\renewcommand\Im{\operatorname{Im}}
\newcommand\tr{\operatorname{tr}}
\newcommand{\SD}{{\mathcal {SD}}}
\newcommand{\SL}{{\mathcal {SL}}}
\newcommand{\slim}{\operatornamewithlimits{s-lim}}

\newcommand\al{\alpha}
\newcommand\be{\beta}
\newcommand\ga{\gamma}
\newcommand\de{\delta}
\newcommand\eps{\varepsilon}
\newcommand\la{\lambda}
\newcommand\si{\sigma}
\newcommand\Si{\Sigma}

\newcommand\fS{{\mathfrak S}}

\newcommand\bC{{\mathbb C}}
\newcommand\bN{{\mathbb N}}
\newcommand\bR{{\mathbb R}}

\newcommand\cH{{\mathcal H}}
\newcommand\cI{{\mathcal I}}
\newcommand\cO{{\mathcal O}}
\newcommand\cP{{\mathcal P}}
\newcommand\cT{\widetilde{T}}

\newcommand\fB{{\mathfrak B}}
\newcommand\fD{{\mathfrak D}}

\title[Inverse spectral problems]{Inverse spectral problems for
    Sturm-Liouville operators with singular potentials${}^{\dag}$}
\author{R.~O.~Hryniv and Ya.~V.~Mykytyuk}

\thanks{${}^{\dag}$The work was partially supported by Ukrainian Foundation
for Basic Research DFFD under grant No.~01.07/00172}

\address{Institute for Applied Problems of Mechanics and Mathematics,
3b~Naukova st., 79601 Lviv, Ukraine and Lviv National University, 1
Universytetska st., 79602 Lviv, Ukraine} \email{rhryniv@iapmm.lviv.ua}

\address{Lviv National University, 1 Universytetska st., 79602 Lviv, Ukraine}
\email{yamykytyuk@yahoo.com}

\subjclass[2000]{Primary 34A55, Secondary 34B24, 34L05, 34L20}
\keywords{Inverse spectral problems, Sturm-Liouville operators, singular
potentials}

\date{November 15, 2002}

\begin{document}

\begin{abstract}
The inverse spectral problem is solved for the class of Sturm-Liouville
operators with singular real-valued potentials from the space
$W^{-1}_2(0,1)$. The potential is recovered via the eigenvalues and the
corresponding norming constants. The reconstruction algorithm is presented
and its stability proved. Also, the set of all possible spectral data is
explicitly described and the isospectral sets are characterized.
\end{abstract}

\maketitle

\section{Introduction}\label{sec:intr}

The main aim of the present paper is to solve the inverse spectral problem
for the class of Sturm-Liouville operators  with singular real-valued
potentials from the space~$W^{-1}_2(0,1)$. Given a real-valued distribution
$q \in W^{-1}_2(0,1)$, we define a Sturm-Liouville operator~$T$ acting in
the Hilbert space $\cH:= L_2(0,1)$ and corresponding to the differential
expression
\begin{equation}\label{eq:S}
    l := - \dd + q
\end{equation}
and, say, the Dirichlet boundary conditions by means of the regularization
method due to \textsc{Savchuk and Shkalikov}~\cite{SS}. Namely, we take a
real-valued $\si \in \cH$ such that $\si' = q$ in the sense of
distributions (thus $\si$ is a distributional primitive of $q$) and put
\begin{equation}\label{eq:Sact}
    T u =T_\si u = l_\si (u) := - (u' -\si u)' - \si u'
\end{equation}
on the domain
\begin{equation}\label{eq:Sdom}
    \fD(T_\si) = \{ u \in W^1_1(0,1) \mid
        u'-\si u \in W^1_1(0,1),\
            l_\si(u) \in \cH ,\ u(0)=u(1)=0\}.
\end{equation}
Observe that, in the sense of distributions, $l_\si(u) = -u'' + q u$ for
all $u\in \fD(T_\si)$. In particular, the operator $T_\si$ does not depend
on the particular choice of the primitive~$\si$ and for \emph{regular}
(i.e., locally summable) potentials it coincides with the standard
Dirichlet Sturm-Liouville operator corresponding to~\eqref{eq:S}.
Also~$T_\si$ depends continuously in the uniform resolvent sense on the
primitive~$\si\in\cH$~\cite{SS} and thus it is a natural Dirichlet
Sturm-Liouville operator related with~\eqref{eq:S} for an arbitrary $q=\si'
\in W^{-1}_2(0,1)$. Note that the class of potentials considered includes,
e.g., the Dirac $\de$-like and Coulomb $1/x$-like potentials that have been
extensively used in quantum mechanics and mathematical
physics~\cite{AGHH,AK}.

It is known~\cite{SS} that for any real-valued $\si\in\cH$ the operator
$T_\si$ defined above is a selfadjoint operator with discrete simple
spectrum $(\la_k^2)$, $k \in \bN$, and that $\la_k$ have the asymptotics
$\la_k = \pi k + \mu_k$ with an $\ell_2$-sequence
$(\mu_k)$~\cite{HM,Sav,SS}. We recall that for regular potentials $q$ the
above asymptotics refines to having $\mu_k = O(1/k)$.

The present paper has arisen as an attempt to answer the following
question: Which sequences $(\la_k^2)$ consisting of real pairwise distinct
numbers and obeying the above-mentioned asymptotics are indeed spectra of
Sturm-Liouville operators with singular potentials from~$W^{-1}_2(0,1)$?
This has led us to the inverse spectral problem for the operators
considered, i.e., to reconstruction of the potential~$q$ based on the
corresponding spectral data.

In the regular situation knowing only the spectrum $(\la_k^2)$ is
insufficient: there are many different potentials $q$ (called {\em
isospectral}) producing Sturm-Liouville operators with the same Dirichlet
spectrum. It was shown by \textsc{P\"oschel and Trubowitz}~\cite{PT} that
the set of all potentials in~$\cH$ with a given admissible spectrum
$(\la_k^2)$ (i.e., real, simple, and obeying the asymptotics $\la_k = \pi k
+ O(1/k)$), is analytically diffeomorphic to the weighted space
$\ell_2(w_n)$ with the weights $w_n = n$. (A rather exceptional situation,
where the spectrum determines the potential uniquely, was pointed out by
\textsc{Am\-bar\-tsu\-myan}~\cite{Am}; namely, he proved that for the
Neumann Sturm-Liouville operator on $(0,1)$ the equalities $\la_k =
\pi(k-1)$ for all $k\in\bN$ imply $q \equiv 0$.)

Therefore to recover the potential $q$ uniquely some additional information
besides the spectrum must be supplied. This can be, e.g., knowledge of the
potential on half of the interval~$(0,1)$~\cite{dRGS,GS,HL,Sakh}, or the
spectrum of a Sturm-Liouville operator given by the same differential
expression but different boundary conditions \cite{Le,Ma}, or three
spectra---one for the whole interval and the others for two halves of
it~\cite{GS1,Pi}. Another kind of information is the squared norms $\al_k$
of properly normalized eigenfunctions (the so-called \emph{norming
constants})~\cite{Le,PT}, and that is exactly the setting we shall treat in
this paper. The aforementioned settings of the inverse spectral problem
will be considered elsewhere.

A complete solution of the inverse spectral problem for a class~$\SL$ of
Sturm-Liouville operators must consist of two parts: (1) an explicit
description of the set $\SD$ of spectral data for the operators in~$\SL$
and (2) development and justification of the method of recovering the
operator in~$\SL$ that corresponds to arbitrary given spectral data
in~$\SD$.

The algorithm of recovering the potential $q$ from the spectral data of a
regular Sturm-Liouville operator based on the transformation operators and
the so-called Gelfand-Levitan-Mar\-chen\-ko equation was developed by
\textsc{Gelfand and Levitan}~\cite{GL} and \textsc{Marchenko}~\cite{Ma1} in
early 1950-ies. In particular, a reconstruction method via the spectrum
$(\la^2_k)$ and the norming constants $(\al_k)$ was suggested in~\cite{GL}
(see also~\cite{Le}), though no precise description was given therein of
the set of all spectral data generated by the considered class of regular
Sturm-Liouville operators. An alternative method for reconstruction of the
potential~$q$ by two spectra was developed by \textsc{Krein}~\cite{K}. A
different approach was suggested later by \textsc{P\"oschel and Trubowitz}
for the class of Sturm-Liouville operators with potentials from~$\cH$
in~\cite{PT}. The authors studied in detail the mapping between the
potentials in~$\cH$ and the spectral data, proved solvability of the
inverse spectral problem and, in particular, completely characterized the
spectral data.

In 1967, \textsc{Zhikov}~\cite{Zh} considered the singular case where the
potential~$q$ is the derivative of a function of bounded variation (i.e.,
where $q$ is a signed measure). The corresponding Sturm-Liouville operator
was defined through the equivalent integral equation, and complete solution
to the inverse spectral problem (in particular, necessary and sufficient
conditions on spectral data) was given. In fact, in~\cite{Zh} the problem
was reduced to recovering $q$ by the corresponding spectral function.

Singularities of different kind were treated by {\sc Rundell and Sacks}
in~\cite{RS}. They considered a Sturm-Liouville operator in impedance form
$Su = \frac1a(au')'$ in the space $L_2((0,1);a)$ with a positive
impedance~$a\in W^1_2[0,1]$. With the help of transformation operators they
found necessary conditions on the spectral data, solved the inverse
spectral problem, and suggested a numerical algorithm for finding a
solution. Observe that $S$ is similar to $Tu=-u''+qu$ with
$q=(\sqrt{a})''/\sqrt{a}$. In particular, for $a\in W^1_2[0,1]$ we get
$q\in W^{-1}_2(0,1)$, though there are many such functions $a$ producing
the same $q\in W^{-1}_2(0,1)$.

\textsc{Coleman and McLaughlin}~\cite{CM} took a different approach to the
inverse spectral problem for Sturm-Liouville operators in impedance form.
Namely, they recast the equation $(au')' + \la au=0$ in the form $u'' +
bu'+\la u=0$ with $b:=a'/a$ and then modified the method of~\cite{PT}
accordingly. In particular, they described the set $\SD$ of the spectral
data and proved unique solvability of the inverse spectral problem for the
case where $b\in L_2(0,1)$.

Note that other types of singularities (e.g.~discontinuous $a$ for
impedance Sturm-Liouville operators $S$, $1/x^\gamma$-like potentials) were
treated by {\sc Anders\-son~\cite{An}, Carlson~\cite{Ca}, Hald~\cite{Ha},
McLaughlin~\cite{ML}, Yurko~\cite{Yu}} a.o.

Here, we generalize the classical approach due to \textsc{Gelfand,
Levitan}, and \textsc{Mar\-chen\-ko} and completely solve the inverse
spectral problem for the class $\SL$ of Sturm-Liouville operators with
singular potentials from~$W^{-1}_2(0,1)$. Namely, we give an explicit
description of the set $\SD$ of spectral data, explain how to recover $q$
from an arbitrary element of $\SD$, and study dependence of $q$ on spectral
data. As a byproduct, we show that the set of potentials in $W^{-1}_2(0,1)$
that are isospectral to a given one is analytically diffeomorphic to the
space $\ell_2$. The main tool of the reconstruction procedure is the
transformation operators for the class of singular Sturm-Liouville
operators in~$\SL$ that were constructed in~\cite{HM}.

Our primary purpose in this paper is to treat the general singular case. We
do not prove here that if the spectral data have better asymptotics, then
the recovered potential is smoother---e.g., if the spectral data formally
correspond to a regular potential from~$\cH$, then the recovered potential
indeed falls into~$\cH$. In a subsequent work we shall justify a more
general claim that the reconstruction algorithm suggested here solves the
inverse spectral problem for the class of Sturm-Liouville operators with
potentials from $W^{\al-1}_2(0,1)$ for every fixed $\al\in[0,1]$.

The organization of the paper is as follows. In the next section we exploit
transformation operators to study the direct spectral problem, i.e., to
give necessary conditions on the set~$\SD$ of spectral data for the
class~$\SL$ of Sturm-Liouville operators considered. In
Section~\ref{sec:cnctn} we establish connection between the spectral data
and the transformation operators and, in particular, derive the so-called
Gelfand-Levitan-Marchenko (GLM) equation. In Section~\ref{sec:GLM} the GLM
equation is proved to possess a unique solution for any element from~$\SD$,
and this solution is shown to be the transformation operator for some
Sturm-Liouville operator from $\SL$ in Section~\ref{sec:inv}. Moreover, the
element from $\SD$ that we have started with turns out to be the spectral
data for the Sturm-Liouville operator found, so that the inverse spectral
problem is completely solved. In Section~\ref{sec:stab} we show stability
of the reconstruction algorithm and characterize isospectral sets, and in
the last section we comment on changes to be made for the case of
Dirichlet-Neumann or third type boundary conditions. Finally, two
appendices contain necessary facts on Riesz bases of sines and cosines in
functional Hilbert spaces and on Hilbert-Schmidt operators respectively.

Throughout the paper, $\si\in\cH$ will denote a real-valued distributional
primitive of the potential $q\in W^{-1}_2(0,1)$, $u^{[1]}$ will stand for
the \emph{quasi-derivative} $u'-\si u$ of a function~$u$, and $\|u\|$ will
denote the $\cH$-norm of $u$.

%
\section{Transformation operators and direct spectral problem}
    \label{sec:direct}
%

In this section we shall solve the direct spectral problem for Dirichlet
Sturm-Liouville operators $T_\si$ with $\si\in \cH$, i.~e., we shall
describe the set $\SD$ of spectral data for $T_\si$---the sequences of
eigenvalues~$(\la^2_k)$ and norming constants~$(\al_k)$ introduced
below---when $\si$ runs over $\cH$. The main tool will be the
transformation operators.

Suppose that $\si \in \cH$ is real-valued; we denote by $\cT_\si$ a
Sturm-Liouville operator $-\dd + \si'$ with the Dirichlet boundary
condition at the point $x=0$. More precisely, $\cT_\si$ acts according to
\[
    \cT_\si u = l_\si(u) = -\bigl(u^{[1]}\bigr)' - \si u^{[1]} - \si^2 u
\]
on the domain
\[
    \fD(\cT_\si) = \{u \in W^1_1[0,1] \mid u^{[1]} \in W^1_1[0,1],
            l_\si(u) \in L_2(0,1), u(0)=0\}.
\]
(We recall that $u^{[1]}:=u'-\si u$ is the quasi-derivative of $u$.)
According to the definition of~$l_\si$, the equation $l_\si(u)=v$ is
understood in the sense that
\[
    \frac{d}{dt}\binom{u}{u^{[1]}} -
    \begin{pmatrix}\si & 1\\-\si^2 & -\si \end{pmatrix}
        \binom{u}{u^{[1]}} = \binom{0}{-v},
\]
so that its solutions enjoy the standard uniqueness properties.

One of the main results established in~\cite{HM} is that the operators
$\cT_\si$ and $\cT_0$ possess the {\em transformation operators}, which
perform their similarity.


\begin{theorem}[\cite{HM}]\label{thm:transf}
 Suppose that $\si \in \cH$; then there exists an integral
 Hilbert-Schmidt operator $K_\si$ of the form
\begin{equation}\label{eq:Kvolt}
    (K_\si u)(x) = \int_0^x k_\si (x,y) u(y)\, dy
\end{equation}
 such that $I + K_\si$ is the transformation operator for
 $\cT_\si$ and $\cT_0$, i.e.,
\begin{equation}\label{eq:simil}
    \cT_\si (I + K_\si) = (I + K_\si) \cT_0.
\end{equation}
If $\si$ is real-valued, then the kernel $k_\si$ of $K_\si$ is real-valued
too. For every $x \in [0,1]$ and $y \in [0,1]$ the functions $k_\si (x,
\cdot )$ and $k_\si (\cdot, y)$ belong to $\cH$ and the mappings $x \mapsto
k_\si (x,\cdot)$ and $y \mapsto k_\si(\cdot,y)$ are continuous in the
$\cH$-norm. Moreover, the operator $K_\si$ with properties
\eqref{eq:Kvolt}--\eqref{eq:simil} is unique.
\end{theorem}


\begin{remark}\label{rem:init}
Some other properties of the transformation operator $K_\si$ established
in~\cite{HM} imply that $I + K_\si$ preserves the initial conditions at
$x=0$ in the sense that for any $u\in W^2_2[0,1]$ and $v:= (I +K_\si)u$ we
have $v(0)=u(0)$ and $v^{[1]}(0) = u'(0)$.
\end{remark}

Denote by $T_\si$ the restriction of $\cT_\si$ by the Dirichlet boundary
condition at $x=1$, i.e., $T_\si := \cT_\si\vert_{\fD_0}$, where $\fD_0 =
\{u \in \fD(\cT_\si) \mid u(1) = 0\}$. It is known~\cite{SS} that $T_\si$
is a bounded below selfadjoint operator with simple discrete spectrum.
Denote by $\la_k^2$, $k \in \bN$, eigenvalues of $T_\si$ in increasing
order. We may assume that all the eigenvalues are positive (and that such
are $\la_k$), as otherwise this situation can be achieved by adding a
suitable constant to $q$.

Let also $u_k$, $k\in\bN$, be the eigenfunction of~$T_\si$ corresponding to
the eigenvalue $\la_k^2$ and normalized by the condition $u_k^{[1]}(0) =
\sqrt{2} \la_k$; then $\al_k:=\|u_k\|^2$ is the corresponding \emph{norming
constant}. Notice that in view of Remark~\ref{rem:init} we have $u_k =
(I+K_\si) v_k$, where $v_k(x) = \sqrt{2}\sin \la_k x$. We also put
$v_{k,0}(x) = \sqrt{2}\sin \pi k x$.

The following statement was proved in~\cite{HM,Sav,SS}, but we sketch its
proof here for the sake of completeness.


\begin{lemma}\label{lem:lak}
 Suppose that $\si\in\cH$ is real-valued and that $\la^2_1 < \la^2_2< \dots$ are
 eigenvalues of $T_\si$ with $\la^2_1>0$. Then $\la_k = \pi k + \mu_k$,
 where the sequence $(\mu_k)_{k=1}^\infty$ belongs to $\ell_2$.
\end{lemma}

\begin{proof}
Since $u_k(x) = v_k(x) + \int_0^x k_\si(x,y) v_k(y)\,dy$, the numbers
$\la_k$ are zeros of the function
\[
    \Phi(\la) := \sin \la + \int_0^1 k_\si(1,y)\sin \la y\, dy.
\]
Recall that $k_\si(1,\cdot) \in \cH$; whence
\(
    \int_0^1 k_\si(1,y)\sin \la y\, dy \to 0
\) as $\la\to\infty$ by the Riemann lemma. Rouch{\'e}'s theorem then gives
$\la_k = \pi k + \mu_k$ with $\mu_k \to 0$ as $k\to\infty$ (see details
in~\cite[Ch.~1.3]{Ma} and \cite{HM}). Therefore  (see
Appendix~\ref{sec:bases}) $\{\sin\la_kx\}_{k=1}^\infty$ is a Riesz basis
of~$\cH$ and the numbers
\[
    \sin\la_k = (-1)^k\sin\mu_k = - \int_0^1 k_\si(1,y)\sin \la_k y\, dy
\]
are the Fourier coefficients of $-k_\si(1,\cdot)\in\cH$ in the biorthogonal
basis. It follows that $(\sin\mu_k)\in \ell_2$, hence $(\mu_k)\in \ell_2$
and the proof is complete.
\end{proof}

The next lemma gives the asymptotics of the norming constants $\al_k$.


\begin{lemma}\label{lem:alk}
Suppose that $\si\in\cH$; then the norming constants $\al_k$ are of the
form $1 + \be_k$, where the sequence $(\be_k)$ belongs to~$\ell_2$.
\end{lemma}

\begin{proof}
Observe that the vectors $u_k$, $k\in\bN$, form a Riesz basis of $\cH$
since so do $v_k$ and $u_k = (I + K_\si)v_k$ with bounded and boundedly
invertible $I + K_\si$. Therefore the $\cH$-norms of $u_k$ are uniformly
bounded, and in view of the relations
\[
    |\be_k| = \bigl|\|u_k\|^2 - \|v_{k,0}\|^2 \bigr| \le
        (1+\|u_k\|)\| u_k - v_{k,0}\|
\]
the lemma will be proved as soon as we show that
\[
    \sum_{k=1}^\infty \| u_k - v_{k,0}\|^2 < \infty.
\]
This inequality follows from the representation
\[
    u_k - v_{k,0} = (I + K_\si)(v_k - v_{k,0})
                + K_\si v_{k,0}
\]
and the fact that both sequences
 $\bigl(\|(I + K_\si)(v_k - v_{k,0})\|\bigr)$
and $\bigl(\|K_\si v_{k,0}\|\bigr)$ belong to $\ell_2$: the former due to
the relation
\[
    v_k(x) - v_{k,0}(x) =
        2\sqrt2\sin \frac{\mu_k}2 \cos(\pi k + \mu_k/2)
        = O(|\mu_k|)
\]
and the inclusion $(\mu_k)\in \ell_2$ (see the previous lemma), and the
latter because $K_\si$ is a Hilbert-Schmidt operator and $(v_{k,0})$ is an
orthonormal basis of $\cH$ (see Appendix~\ref{sec:HS}). The proof is
complete.
\end{proof}

Denote by $\SD$ the set of all pairs of sequences
$\{(\la^2_k)_{k=1}^\infty,(\al_k)_{k=1}^\infty\}$ satisfying the following
two conditions:
\begin{itemize}
\item[(A1)] $\la_k$ are all positive, strictly increase with $k$, and obey
    the asymptotic relation $\la_k = \pi k + \mu_k$ with some $\ell_2$-sequence
        $(\mu_k)_{k=1}^\infty$;
\item[(A2)] $\al_k = 1 + \be_k > 0$ with some $\ell_2$-sequence
    $(\be_k)_{k=1}^\infty$.
\end{itemize}
Also $\SL$ will stand for the set of all positive Dirichlet Sturm-Liouville
operators $T_\si$ with $\si\in\cH$ real-valued. Lemmata~\ref{lem:lak} and
\ref{lem:alk} demonstrate that the spectral data for the operators in $\SL$
belong to $\SD$. Our next task is to show that, conversely, any element of
$\SD$ is spectral data for some operator in $\SL$.

%
\section{Connection between spectral data and transformation operators}%
    \label{sec:cnctn}
%

In this section we shall derive a relation between the spectral data for a
Dirichlet Sturm-Liouville operator $T_\si$ with $\si\in\cH$ and the
transformation operator $K_\si$. This relation will be used in the next
section to find $K_\si$ given the spectral data and thus to solve the
inverse spectral problem.

Suppose therefore that $\{(\la_k^2),(\al_k)\}\in\SD$ is the spectral data
for a Sturm-Liouville operator $T_\si$ with $\si\in \cH$. We recall that
$\la_k = \pi k + \mu_k>0$ and $\al_k = 1 + \be_k>0$ with some
$\ell_2$-sequences $(\mu_k)$ and $(\be_k)$. Put
\[
   U := \slim_{N\to\infty} \sum_{k=1}^N \frac1{\al_k} (\cdot, v_k)v_k,
\]
where $v_k(x) = \sqrt{2}\sin \la_k x$ and $\slim$ stands for the limit in
the strong operator topology of $\cH$. Since by Proposition~\ref{pro:RB}
the system $(v_k)_{k=1}^\infty$ forms a Riesz basis of $\cH$ and
$\inf_{k\in\bN} 1/\al_k > 0$, the operator $U$ is bounded, selfadjoint, and
uniformly positive.


\begin{lemma}\label{lem:u}
For all $j,k\in\bN$, the following relation holds: $(U^{-1}v_j,v_k) = \al_k
\de_{jk}$, where $\de_{jk}$ is the Kronecker delta.
\end{lemma}

\begin{proof}
Note that
\[
    U^{-1} = \slim_{N\to\infty} \sum_{k=1}^N {\al_k} ( \cdot, w_k)w_k,
\]
where $(w_k)$ is a basis biorthogonal to $(v_k)$ (see~\cite[Ch.~VI]{GK1}
and also Appendix~\ref{sec:bases}). Therefore
\[
    (U^{-1}v_j,v_k) = \slim_{N\to\infty}
        \sum_{l=1}^N {\al_l} ( v_j , w_l)(w_l, v_k)
    = \al_l \de_{jl}\de_{lk} = \al_k \de_{jk},
\]
and the proof is complete.
\end{proof}


\begin{lemma}\label{lem:F}
$F:=U-I$ is an integral operator of the Hilbert-Schmidt class with kernel
\begin{equation}\label{eq:F}
    f(x,y) = \phi (x+y) - \phi (x-y),
\end{equation}
where
\begin{equation}\label{eq:phi}
   \phi(s) = \sum_{k\in\bN}
        \Bigl( \cos \pi ks - \frac1\al_k \cos\la_ks\Bigr)
\end{equation}
is an $L_2(0,2)$-function.
\end{lemma}

\begin{proof}
Recall that we have denoted by $v_{k,0}$, $k\in\bN$, the function
$\sqrt{2}\sin\pi kx$. Since
 \(
    I = \slim_{N\to\infty} \sum_{k=1}^N (\cdot,v_{k,0}) v_{k,0},
 \)
we have
\[
    U-I = \slim_{N\to\infty} \sum_{k=1}^N
    \Bigl( \frac1{\al_k} ( \cdot, v_k)v_k
        - (\cdot,v_{k,0}) v_{k,0} \Bigr).
\]
Observe that
 \(
  \dfrac1{\al_k} ( \cdot, v_k)v_k
        - (\cdot,v_{k,0}) v_{k,0}
 \)
is an integral operator with kernel
\begin{multline*}
    f_k(x,y) = \frac2{\al_k} \sin\la_kx\,\sin\la_ky -
        2\sin\pi kx \sin \pi ky \\
        = \cos \pi k(x+y) - \frac1{\al_k} \cos\la_k(x+y)
        - \cos \pi k(x-y) + \frac1{\al_k} \cos\la_k(x-y).
\end{multline*}
Letting
\[
    \phi_N :=  \sum_{k=1}^N
        \Bigl( \cos \pi ks - \frac1\al_k \cos\la_ks\Bigr),
\]
we see that
\[
    F_N := \sum_{k=1}^N
    \Bigl( \frac1{\al_k} ( \cdot, v_k)v_k
        - (\cdot,v_{k,0}) v_{k,0} \Bigr)
\]
is an integral operator with kernel $f_N(x,y):=\phi_N (x+y) - \phi_N(x-y)$.

Put $\ti\be_k:=\be_k/\al_k$; then $1/\al_k = 1 - \ti\be_k$ with
$(\ti\be_k)_{k=1}^\infty \in \ell_2$ and
\[
    \cos \pi ks - \frac1{\al_k} \cos \la_ks =
        2 \sin (\mu_k s/2) \sin [(\pi k + \mu_k/2)s]
        + \ti\be_k \cos\la_ks.
\]
Recall (see Corollary~\ref{cor:RB}) that both
 $\bigl(\sin[(\pi k + \mu_k/2)s]\bigr)_{k=1}^\infty$ and
$(\cos\la_ks)_{k=1}^\infty$ constitute Riesz bases of their closed linear
spans in $L_2(0,2)$. Since, moreover,
\[
 2 \sin (\mu_k s/2) = s \mu_k + O(|\mu_k|^3)
\]
as $k\to\infty$, the series for $\phi$ is convergent in $L_2(0,2)$ and
hence $\phi_N\to\phi\in L_2(0,2)$. Thus $f_N$ converge in the
$L_2\bigl((0,1)\times(0,1)\bigr)$-norm to $f$ given by~\eqref{eq:F}. This
implies that $F_N$ converge in the Hilbert-Schmidt norm to the integral
operator $F$ with kernel $f$, and the lemma is proved.
\end{proof}


\begin{lemma}\label{lem:id}
$(I+K_\si) (I + F) (I+ K_\si^*) = I$.
\end{lemma}

\begin{proof}
Since $(u_k/\sqrt{\al_k})_{k=1}^\infty$ is an orthonormal basis of $\cH$
and $u_k = (I + K_\si) v_k$, we have
\begin{multline*}
    I = \slim_{N\to\infty} \sum_{k=1}^N \frac1{\al_k} (\cdot,u_k) u_k
      = \slim_{N\to\infty} \sum_{k=1}^N \frac1{\al_k} (\,\cdot\,,(I+ K_\si)v_k)
                    (I+K_\si)v_k \\
    = (I + K_\si) \Bigl[
        \slim_{N\to\infty} \sum_{k=1}^N \frac1{\al_k} (\cdot,v_k) v_k
    \Bigr] (I+K_\si^*)
    = (I+K_\si) (I + F) (I+K_\si^*).
\end{multline*}
The lemma is proved.
\end{proof}


\begin{corollary}\label{cor:GLM}
The kernels $k_\si(x,y)$ and $f(x,y)$ of the operators $K_\si$ and $F$
satisfy the following Gelfand-Levitan-Marchenko (GLM) equation for a.e.
$x,y$ with $x>y$:
\begin{equation}\label{eq:GLM}
    k_\si(x,y) + f(x,y) + \int_0^x k_\si(x,s) f(s,y)\,ds = 0.
\end{equation}
\end{corollary}

\begin{proof}
Observe that $K_\si^*$ is an integral Volterra operator with upper-diagonal
kernel $k_\si^*(x,y) = \ov{k_\si(y,x)}$, i.e.,
\[
    (K_\si^* u)(x) = \int_x^1 {k_\si^* (x,y)} u(y)\,dy
\]
and $k_\si^*(x,y)=0$ for $x>y$. Therefore the operator $I + K_\si^*$ is
invertible and its inverse can be written in the form $I + \ti K$, where
$\ti K$ is an integral Volterra operator with upper-diagonal kernel. By
Lemma~\ref{lem:id},
\[
    (I+K_\si)(I+F) = (I + K_\si^*)^{-1} = I + \ti K;
\]
spelling out this equality in terms of the kernels $k_\si$ and $f$ for
$x>y$, we arrive at the GLM equation~\eqref{eq:GLM}.
\end{proof}

%
\section{Solution of the GLM equation}\label{sec:GLM}
%

In the previous section, we showed that the spectral data
$\{(\la_k^2),(\al_k)\}$ for the Sturm-Liouville operator $T_\si$ are
connected with the transformation operator $K_\si$ for the pair $\cT_\si$
and $\cT_0$ through the GLM equation~\eqref{eq:GLM}. In this section we
shall prove that the GLM equation is uniquely soluble for $K_\si$.

Recall (see Appendix~\ref{sec:HS}) that the ideal $\fS_2$ of
Hilbert-Schmidt operators consists of integral operators with square
summable kernels and that the scalar product $\lan X,Y\ran_2 := \tr(X Y^*)$
introduces a Hilbert space topology on $\fS_2$. Denote by $\fS_2^+$ the
subspace of $\fS_2$ consisting of all Hilbert-Schmidt operators with
lower-triangular kernels. In other words, $T \in \fS_2$ belongs to
$\fS_2^+$ if the kernel $t$ of~$T$ satisfies $t(x,y) = 0$ for $0 \le x < y
\le 1$. For an arbitrary $T\in\fS_2$ with kernel $t$ the cut-off $t^+$ of
$t$ given by
\[
    t^+(x,y) = \left\{ \begin{array}{ll}
        t(x,y)& \quad \mbox {for $x \ge y$}\\
              0 & \quad \mbox {for $x < y$} \end{array}
        \right.
\]
generates an operator $T^+\in\fS_2^+$, and the corresponding mapping
$\cP^+: T \mapsto T^+$  turns out to be an orthoprojector in $\fS_2$ onto
$\fS_2^+$, i.e. $(\cP^+)^2 = \cP^+$ and
 $\lan \cP^+ X, Y\ran_2 = \lan X,\cP^+ Y\ran_2$ for any $X,Y\in \fS_2$; see
 details in~\cite[Ch.~I.10]{GK2}.

With these notations, the GLM equation~\eqref{eq:GLM} can be recast as
\begin{equation}\label{eq:KGLM}
    K + \cP^+ F + \cP^+(K F) = 0
\end{equation}
or
\[
    (\cI + \cP^+_F) K = - \cP^+ F,
\]
where $\cP^+_X$ is the linear operator in $\fS_2$ defined by $\cP^+_X Y =
\cP^+ (YX)$ and $\cI$ is the identity operator in $\fS_2$. Therefore
solvability of the GLM equation is closely connected with the properties of
the operator~$\cP^+_F$.


\begin{lemma}\label{lem:Pcnts}
The operator $\cP^+_X$ depends continuously on $X \in \fB(\cH)$.
\end{lemma}

\begin{proof}
This is a straightforward consequence of the fact that $\cP^+_X$ depends
linearly on~$X$ and that 
\[
    \|\cP^+_X Y\|_{\fS_2} \le \|YX\|_{\fS_2} \le \|X\| \|Y\|_{\fS_2},
\]
see~\cite[Ch.~3]{GK1}.
\end{proof}


\begin{lemma}\label{lem:Ppos}
Suppose that $I+X$ is a uniformly positive operator in $\cH$. Then $\cI +
\cP^+_X$ is uniformly positive in $\fS_2^+$.
\end{lemma}

\begin{proof}
Given any $Y\in\fS_2^+$, we observe that
\[
    \lan Y + \cP^+_X Y, Y\ran_2 = \lan Y,Y\ran_2 + \lan YX,Y \ran_2
        = \tr \bigl(Y(I+X)Y^*\bigr).
\]
If $\eps>0$ is chosen so that $I+X \ge \eps I$, then $Y(I+X)Y^* \ge \eps
YY^*$ and by monotonicity of the trace we get
\[
    \lan (\cI+\cP^+_X) Y,Y \ran_2 \ge \eps \lan Y, Y \ran_2.
\]
The lemma is proved.
\end{proof}


\begin{lemma}\label{lem:solveGLM}
Suppose that sequences $(\la_k^2)_{k=1}^\infty$ and $(\al_k)_{k=1}^\infty$
satisfy conditions (A1) and (A2) and that $F$ is an integral operator with
kernel $f$ of~\eqref{eq:F}. Then equation~\eqref{eq:KGLM} has a unique
solution $K$, which belongs to $\fS_2^+$.
\end{lemma}

\begin{proof}
Under assumptions (A1) and (A2) the operator $F$ is of Hilbert-Schmidt
class and
\[
    I + F = \slim_{N\to\infty}
        \sum_{k=1}^N \frac1{\al_k} (\,\cdot\, , v_k) v_k
\]
(cf.~the proof of Lemma~\ref{lem:F}). Moreover, $U:=I+F$ is uniformly
positive as noted at the beginning of Section~\ref{sec:cnctn}. Therefore
the operator $\cI+\cP^+_F$ is uniformly positive in $\fS_2^+$ by
Lemma~\ref{lem:Ppos} and
\[
    K := - (\cI + \cP^+_F)^{-1} \cP^+ F \in \fS_2^+
\]
is a unique solution of~\eqref{eq:KGLM}. The lemma is proved.
\end{proof}


\begin{remark}
It is well known that a necessary and sufficient condition for the GLM
equation to be soluble is that, for any $a\in [0,1]$, the operator
\[
    u \mapsto u(x) + \int_0^a f(x,y)u(y)\,dy
\]
is boundedly invertible in $\cH$. For a symmetric kernel $f$ this condition
is easily seen to be equivalent to uniform positivity of the operator
$I+F$.

Also, the GLM equation is uniquely soluble if and only if $I+F$ can be
factorized as $(I + X^+)(I + X^-)$ with some $X^+ \in \fS_2^+$ and $X^- \in
\fS_2^- := \fS_2 \ominus \fS_2^+$, see~\cite[Ch.~IV]{GK2}.
\end{remark}

%
    \section{The inverse spectral problem}\label{sec:inv}
%

In the previous two sections we showed that the operator $F$ constructed
via the spectral data $\{(\la_k^2), (\al_k)\}$ of a Dirichlet
Sturm-Liouville operator $T_\si$, $\si\in\cH$, is connected with the
transformation operator $I + K_\si$ for $\cT_\si$ and $\cT_0$ through the
GLM equation~\eqref{eq:GLM} or \eqref{eq:KGLM} and then proved that the GLM
equation is uniquely soluble for~$K_\si$. Suppose that instead we have
started with an arbitrary element $\{(\la_k^2), (\al_k)\} \in \SD$,
constructed $F$ in the same manner, and found a solution $K$ of the GLM
equation; is then $I+K$ a transformation operator for some $\cT_\si$,
$\si\in\cH$, and $\cT_0$?

In this section we give an affirmative answer to this question and find an
explicit formula for the corresponding $\si$. Moreover, we show that the
element $\{(\la_k^2), (\al_k)\}$ we have started with is indeed the
spectral data for the Sturm-Liouville operator~$T_\si$ with the $\si$
found. This completes the solution of the inverse spectral problem for the
considered class of Sturm-Liouville operators.


\begin{theorem}\label{thm:findK}
Suppose that $\{(\la_k^2), (\al_k)\}$ is an arbitrary element of $\SD$ and
that $F$ is an integral operator with kernel $f$ of~\eqref{eq:F}. Let also
$K$ be a (unique) solution of the GLM equation~\eqref{eq:KGLM}. Then there
exists a (unique up to an additive constant) function $\si \in \cH$ such
that $I+K$ coincides with the transformation operator $I + K_\si$ for the
pair of Sturm-Liouville operators $\cT_\si$ and $\cT_0$.
\end{theorem}

\begin{proof}
We shall approximate $F$ in the $\fS_2$-norm by a sequence
$(F_n)_{n=1}^\infty$ of Hilbert-Schmidt operators with smooth (say,
infinitely differentiable) kernels $f_n$ so that the following holds:
\begin{itemize}
\item[(a)] the solutions $K_n$ of~\eqref{eq:KGLM} for $F$ replaced with
        $F_n$ converge to $K$ in $\fS_2$ as $n\to\infty$;
\item[(b)] for each $n\in\bN$ there exists $\si_n\in \cH$ such that
    $I+ K_n$ is a transformation operator for the pair $\cT_{\si_n}$
    and $\cT_0$;
\item[(c)] $(\si_n)_{n=1}^\infty$ is a Cauchy sequence in $\cH$.
\end{itemize}
Put $\si:=\lim \si_n$; then by~\cite{HM} the operators $I+K_n$ converge
in~$\cH$ to an operator $I + K_\si$, which is the transformation operator
for the pair $\cT_\si$ and $\cT_0$. Thus $K = K_\si$ yielding the result.
The uniqueness of $\si$ up to an additive constant is obvious.

The details are as follows. We put
\[
    f_n(x,y) := \phi_n(x+y) - \phi_n(x-y),
\]
where
\[
    \phi_n(s) := \sum_{k=1}^n
        \Bigl( \cos \pi ks - \frac1{\al_k} \cos\la_ks\Bigr).
\]
We recall (see the proof of Lemma~\ref{lem:F}) that $\phi_n \to \phi$ in
$\cH$ as $n \to \infty$. This implies that $f_n$ converge to $f$ in
$L_2\bigl((0,1)\times(0,1)\bigr)$, i.e., that $F_n$ converge to $F$ in
$\fS_2$ as $n \to \infty$. The corresponding operators $I + F_n$ are
uniformly positive so that the equations
    $K_n + \cP^+ F_n + \cP^+(K_nF_n) = 0$
have unique solutions
    $K_n = - (\cI + \cP^+_{F_n})^{-1} \cP^+ F_n\in\fS_2^+$.
Combining all these statements, we find that
\[
    \|(\cI + \cP^+_{F_n})^{-1} - (\cI + \cP^+_{F})^{-1}\|_{\fB(\fS_2)}
        \to 0
\]
and
\begin{multline*}
   \|K_n - K\|_{\fS_2} =
    \|(\cI + \cP^+_{F_n})^{-1} \cP^+F_n -
        (\cI + \cP^+_{F})^{-1}\cP^+ F \|_{\fS_2}\\
    \le \|(\cI + \cP^+_{F_n})^{-1} - (\cI + \cP^+_{F})^{-1}\|_{\fB(\fS_2)}
        \|\cP^+F_n\|_{\fS_2} \\
    + \|(\cI + \cP^+_{F})^{-1}\|_{\fB(\fS_2)}
        \|\cP^+F_n - \cP^+F\|_{\fS_2} \to 0
\end{multline*}
as $n \to \infty$ thus establishing (a).

Note that the kernel $k_n$ of $K_n$ solves the GLM equation
\[
    k_n(x,y) + f_n(x,y) + \int_0^x k_n(x,s) f_n(s,y)\,ds =0,
        \qquad x>y.
\]
Since $f_n$ is smooth, the classical result \cite[Ch.~II]{Le},
\cite[Ch.~2.3]{Ma} states that $k_n$ is at least once continuously
differentiable and that $I + K_n$ is the transformation operator for the
Sturm-Liouville operators $\cT_{\si_n}$ and $\cT_0$, where $\si_n$ is a
primitive of $q_n(x) := 2 \frac{d}{dx}k_n(x,x)$. The suitable choice of the
primitive (suggested by the GLM equation above) is
\begin{equation}\label{eq:sigman}
    \si_n(x) := 2k_n(x,x) + 2 \phi_n(0) =
        - 2 \phi_n(2x) - 2\int_0^x k_n(x,s) f_n(s,x)\,ds,
\end{equation}
and (b) is fulfilled.

To prove (c), we observe that
\begin{align*}
    \|\si_n - \si_m\|^2 &\le 12\int_0^1|\phi_n(2x) - \phi_m(2x)|^2\,dx \\
        &+ 12\int_0^1 dx \Bigl| \int_0^x \bigl(k_n(x,y) - k_m(x,y)\bigr)
                    f_n(y,x)\,dy\Bigr|^2 \\
        &+ 12\int_0^1 dx \Bigl| \int_0^x k_m(x,y)
                \bigl(f_n(y,x) - f_m(y,x)\bigr)\,dy\Bigr|^2.
\end{align*}
Since $f_n(x,y) = \phi_n(x+y) - \phi_n(x-y)$ and $\phi_n$ form a Cauchy
sequence in $L_2(0,2)$, we find that
\begin{equation}\label{eq:Fcnt}
\begin{aligned}
  \int_0^x |f_n(y,x)|^2\,dy &\le
     2\int_0^x |\phi_n(y+x)|^2\,dy +  2\int_0^x |\phi_n(y-x)|^2\,dy \\
     &= 2 \int_0^{2x} |\phi_n(s)|^2\,ds \le 2 \int_0^2 |\phi_n(s)|^2\,ds
     \le C
\end{aligned}
\end{equation}
for some positive $C$ independent of $n\in\bN$ and $x\in[0,1]$, and hence
\begin{align*}
    \int_0^1 dx \Bigl(
        \int_0^x\bigl(k_n(x,s) &- k_m(x,s)\bigr)f_n(s,x)\,ds
            \Bigr)^2 \\
    &\le \int_0^1 dx \Bigl(
        \int_0^x\bigl|k_n(x,s) - k_m(x,s)\bigr|^2\,ds
        \int_0^x |f_n(s,x)|^2\,ds \Bigr) \\
    &\le C \int_0^1 \int_0^1 |k_n(x,s) - k_m(x,s)\bigr|^2\,ds\,dx\\
    &= C \|K_n - K_m\|^2_{\fS_2} \to 0
\end{align*}
as $n,m\to\infty$. In the same manner it can be shown that
\[
    \int_0^x |f_n(s,x) - f_m(s,x)|^2\,ds
        \le 2 \int_0^2 |\phi_n(s) - \phi_m(s)|^2\,ds
        =2\|\phi_n - \phi_m\|^2_{L_2(0,2)}
\]
and
\[
   \int_0^1 dx\Bigl(\int_0^x k_m(x,s) \bigl(f_n(s,x) - f_m(s,x)\bigr) \,ds\Bigr)^2
    \le 2 \|K_m\|^2_{\fS_2} \|\phi_n - \phi_m\|^2_{L_2(0,2)} \to 0
\]
as $n,m \to \infty$.
Combining the above relations, we conclude that $(\si_n)$ is a Cauchy
sequence in $\cH$. Therefore (c) is satisfied and the proof is complete.
\end{proof}

\begin{remark}\label{rem:si}
The arguments used to prove (c) above also justify passage to the limit in
the $\cH$-sense in~\eqref{eq:sigman}, and this results in the equality
\begin{equation}\label{eq:sigma}
    \si(x) = -2\phi(2x) - 2 \int_0^x k(x,s) f(s,x)\, ds,
\end{equation}
where $k$ is the kernel of $K$. If $k$ and $\phi$ are smooth, then the GLM
equation implies that \(
    \si(x) = 2 k(x,x) - 2 \phi(0),
\)
thus yielding the classical relation
\[
    q(x) = 2 \frac{d}{dx} k(x,x)
\]
for the potential $q$.
\end{remark}


To complete the arguments, we have to show that the spectral data for
$T_\si$ with $\si\in\cH$ just found coincide with the data $\{(\la_k^2),
(\al_k)\}\in\SD$ that we have started with.


\begin{theorem}\label{thm:spdt}
 $\{(\la_k^2), (\al_k)\}$ are the spectral data for $T_\si$.
\end{theorem}

\begin{proof}
Put $w_k:= (I + K) v_k$; then $w_k(0)=0$, $w^{[1]}_k(0)=\sqrt2\la_k$, and,
due to similarity of $\cT_\si$ and $\cT_0$, $\cT_\si w_k = \la_k^2 w_k$.
The system $(w_k)_{k=1}^\infty$ is complete in~$\cH$ since such is
$(v_k)_{k=1}^\infty$ and $I+K$ is a homeomorphism. Also,
Lemmata~\ref{lem:u} and \ref{lem:id} imply the orthogonality relation
\[
    (w_j,w_k) = \bigl((I+K^*)(I + K) v_j,v_k \bigr)
    = (U^{-1} v_j, v_k) = \al_k \de_{jk};
\]
in particular, $\al_k = \|w_k\|^2$.

It remains to show that $w_k$ are eigenfunctions of $T_\si$, i.e., that
$w_k(1) = 0$. Observe that
\begin{equation}\label{eq:lagr}
    0 = (\cT_\si w_j,w_k) - (w_j, \cT_\si w_k)
     = - w_j^{[1]}(1) \ov{w_k(1)} +
        w_j(1) \ov{w_k^{[1]}(1)}.
\end{equation}
If $w_j(1) = 0$ for some $j\in\bN$, then $w_j^{[1]}(1)\ne0$ due to
uniqueness of solutions of the equation $l_\si(y) = \la_j^2 y$, and the
above relation shows that $w_k(1) = 0$ for all $k\in\bN$ as required.

Otherwise $w_j(1)\ne0$ for all $j\in\bN$, and equation~\eqref{eq:lagr}
implies that there exists a constant $h \in \bR$ such that
 \(
    w_j^{[1]}(1)/w_j(1) = h
 \)
for all $j \in \bN$. Thus $\la^2_k$ are eigenvalues of the Sturm-Liouville
operator $T_{\si,h}$, which is the restriction of $\cT_\si$ by the boundary
condition $y^{[1]}(1) = h y(1)$. However, the eigenvalues~$\la^2_k$ of
$T_{\si,h}$ are known to obey the asymptotics $\la_k = \pi(k-1/2) + o(1)$
as $k\to \infty$ (see Section~\ref{sec:third}), and this contradiction
eliminates the possibility that $w_j(1)$ do not vanish.

Thus we have proved that every function $w_k$ satisfies the Dirichlet
condition at $x=1$ and so is an eigenfunction of the Sturm-Liouville
operator~$T_\si$ corresponding to the eigenvalue~$\la_k^2$. Moreover,
$T_\si$ has no other eigenvalues since the system~$(w_k)_{k=1}^\infty$ is
already complete in~$\cH$. Henceforth the element $\{(\la_k^2),
(\al_k)\}\in\SD$ is indeed the spectral data for the Sturm-Liouville
operator~$T_\si$ with $\si$ of~\eqref{eq:sigma}, and the theorem is proved.
\end{proof}

To sum up, we have established the following:

\begin{corollary}\label{rem:isp}
Two sequences $(\la_k^2)$ and $(\al_k)$ are the spectral data for a
(positive) Sturm-Liouville operator $T_\si$ with potential $q = \si'$ from
the space $W^{-1}_2(0,1)$ if and only if assumptions (A1) and (A2) are
satisfied.

The corresponding operator $T_\si$ is uniquely recovered from the spectral
data through formula~\eqref{eq:sigma}, in which the functions $\phi$ and
$f$ are given by~\eqref{eq:phi} and \eqref{eq:F} respectively, and $k$ is
the kernel of the solution $K$ of the GLM equation~\eqref{eq:KGLM}.
\end{corollary}

%
\section{Stability and isospectral sets}\label{sec:stab}
%

In this section, we would like to study the correspondence between the
Dirichlet Sturm-Liouville operators $T_\si$ with real-valued $\si\in\cH$
and their spectral data in more detail. Namely, we shall show that the
potential $q=\si'\in W^{-1}_2(0,1)$ (and thus the operator $T_\si\in\SL$)
depends continuously on the spectral data $\{(\la_k^2),(\al_k)\}\in\SD$ and
that the isospectral sets are analytically diffeomorphic to the Hilbert
space $\ell_2$. Similar results were established in~\cite{ML,PT} for the
regular case $q\in \cH$, and in~\cite{An,CM} for impedance Sturm-Liouville
operators.

First we have to introduce the topology on the set $\SD$. Recall that by
definition any element $\{(\la_k^2),(\al_k)\}\in\SD$ is uniquely determined
by two $\ell_2$-sequences $(\mu_k)$ and $(\be_k)$ through the relations
$\la_k = \pi k + \mu_k$ and $\al_k = 1 + \be_k$. Therefore we can identify
$\SD$ with an open subset of the space $\ell_2\times \ell_2$ in the
standard coordinate system $\bigl((\mu_k),(\be_k)\bigr)$. In this way $\SD$
becomes a subset of a Hilbert space and inherits the topology of that
space.

We shall study the correspondence between spectral data
$\{(\la_k^2),(\al_k)\}\in\SD$ and operators $T_\si\in\SL$ through the chain
\[
    \{(\la_k^2),(\al_k)\} \mapsto \phi \mapsto K \mapsto \si \mapsto T_\si,
\]
in which $\phi\in L_2(0,2)$ is the function of~\eqref{eq:phi} and
$K\in\fS_2$ is the operator of~\eqref{eq:KGLM}.


\begin{lemma}\label{lem:phicont}
The function
\[
   \phi(s) = \sum_{k\in\bN}
        \Bigl( \cos \pi ks - \frac1\al_k \cos\la_ks\Bigr)
\]
depends continuously in $L_2(0,2)$ on the spectral data
$\{(\la_k^2),(\al_k)\}\in\SD$.
\end{lemma}

\begin{proof}
We can rewrite the function $\phi$ in terms of the $\ell_2$-sequences
$(\mu_k)$ and $(\ti\be_k)$ with $\ti\be_k:= 1-1/\al_k = \be_k/\al_k$ as
follows:
\begin{align*}
    \phi(s) &= \sum_{k\in\bN} \bigl( \cos\pi ks - \cos(\pi ks + \mu_ks)
        + \ti\be_k \cos(\pi ks + \mu_ks) \bigr)\\
        &= \sum_{k\in\bN} 2\sin(\mu_ks/2)\sin[(\pi k+ \mu_k/2)s]
            + \sum_{k\in\bN} \ti\be_k \cos(\pi ks + \mu_ks).
\end{align*}
Therefore it suffices to prove that the mappings from $\ell_2\times \ell_2$
into $L_2(0,2)$ given by
\[
    \{(\nu_k),(\ga_k)\} \mapsto \psi_1(s):=
            \sum_{k\in\bN} \ga_k  \cos(\pi ks + \nu_ks)
\]
and
\[
    \{(\nu_k),(\ga_k)\} \mapsto \psi_2(s):=
            \sum_{k\in\bN} \ga_k  \sin(\pi ks + \nu_ks)
\]
are continuous.

Suppose that $\{(\hat\nu_k),(\hat\ga_k)\}$ is another element
of~$\ell_2\times \ell_2$ giving rise to a function
    $\hat \psi_1 := \sum_{k\in\bN} \hat\ga_k
        \cos(\pi ks + \hat\nu_ks) \in L_2(0,2)$;
then
\begin{equation}\label{eq:differ}
\begin{aligned}
    \|\psi_1 - \hat \psi_1\|_{L_2(0,2)}
        &\le \|\sum (\ga_k - \hat\ga_k) \cos(\pi ks+\nu_ks)\|_{L_2(0,2)} \\
        & + \|\sum \hat\ga_k[\cos(\pi ks+\nu_ks) -
                \cos(\pi ks+\hat\nu_ks)]\|_{L_2(0,2)}.
\end{aligned}
\end{equation}
We may assume that the sequence $(\pi k + \nu_k)$ strictly increases; then
$\bigl( \cos(\pi ks+\nu_ks) \bigr)_{k=1}^\infty$ is a Riesz basic sequence
in~$L_2(0,2)$, and the first summand above is bounded by
 \(
   C \|(\ga_k)-(\hat\ga_k)\|_{\ell_2}
 \)
with $C=C\bigl((\nu_k)\bigr)$ being the corresponding Riesz constant (see
Appendix~\ref{sec:bases}). In view of the inequality
\[
    |\cos(\pi ks + \nu_ks) - \cos(\pi ks + \hat\nu_ks)|
                \le  |\nu_k - \hat\nu_k|, \qquad s\in(0,2),
\]
the second summand in~\eqref{eq:differ} is bounded by
 \(
        \|(\hat \ga_k)\|_{\ell_2} \|(\nu_k)-(\hat\nu_k)\|_{\ell_2},
 \)
and the statement for $\psi_1$ is proved. Continuity of the second mapping
is proved analogously.
\end{proof}


\begin{lemma}\label{lem:Kcont}
The solution $K\in\fS_2$ of the GLM equation~\eqref{eq:KGLM} depends
locally analytically on $\phi\in L_2(0,2)$.
\end{lemma}

\begin{proof}
The operator $F$ depends linearly on $\phi\in L_2(0,2)$ and, in view of
inequality~\eqref{eq:Fcnt},
 $\|F\|\le\|F\|_{\fS_2}\le 2\|\phi\|_{L_2(0,2)}$,
so that $F$ is an analytic function of $\phi$ (see~\cite[Ch.~2]{D} or
\cite[App.~A]{PT} on analytic mappings of Banach spaces). Next, by
Lemma~\ref{lem:Pcnts} the operator~$\cP_F^+\in\fB(\fS_2)$ is continuous in
$F\in\fS_2$ (and thus analytic in view of linearity), while the inverse
function $(\cI - \cP_F^+)^{-1}$ is locally analytic in $\cP_F^+$. Combining
these statements and recalling the formula
\[
    K = - (\cI + \cP^+_F)^{-1} \cP^+ F,
\]
we easily derive the claim.
\end{proof}


\begin{lemma}\label{lem:sicont}
The function $\si$ of~\eqref{eq:sigma} depends locally analytically
in~$\cH$ on $\phi$ in $L_2(0,2)$.
\end{lemma}

\begin{proof}
By definition
 \(
    \si(x) = -2\phi(2x) - 2 \int_0^x k(x,s) f(s,x)\, ds,
 \)
where $k$ is the kernel of $K$ and $f(s,x)=\phi(s+x)-\phi(s-x)$. The second
summand above is a continuous bilinear function of~$K$ and~$\phi$ (see the
proof of Theorem~\ref{thm:findK}). Therefore~$\si$ is a jointly analytic
function of~$\phi$ and~$K$~\cite[Ch.~2]{D}, and in view of
Lemma~\ref{lem:Kcont} the result follows.
\end{proof}


We denote by $\Si^+$ the set of all real-valued $\si\in\cH$, for which the
operators $T_\si$ are positive. Observe that the operators $T_\si$, $K_\si$
and the spectral data do not change if $\si$ is replaced by $\si+c$ with
any real $c$, so that these objects depend in fact on an equivalence class
$\hat\si$ in $\Si^+/\bR$ (i.e., on the common derivative $q=\si'$ for
$\si\in\hat\si$) rather than on $\si\in\Si^+$. It is proved in~\cite{HM,SS}
that this dependence is continuous for $T_\si$ and $K_\si$; we shall show
next that also the spectral data depend continuously on $\hat\si \in
\Si^+/\bR$.

\begin{lemma}[{cf.~\cite[Lemma~4.3]{An}}]
The mapping
 \(
    \Si^+ /\bR \ni \hat\si \mapsto (\la_k^2)
 \)
is continuous.
\end{lemma}

\begin{proof}
We recall that $\la_k=\pi k+\mu_k$ for some $\ell_2$-sequence $(\mu_k)$ and
that the continuity of~$(\la_k^2)$ is understood as continuity of~$(\mu_k)$
in the~$\ell_2$-topology. Also, $\la_k$ are simple zeros of the function
\[
    \Phi(\la)=\Phi(\la,\si):= \sin\la + \int_0^1 k_\si(1,y)\sin\la y\, dy,
\]
where $k_\si$ is the kernel of the corresponding transformation
operator~$K_\si$.

It can be shown (see also~\cite[Ch.~1.3]{Ma}) that
\begin{align*}
    |\Phi(\la) - \sin\la| &= \Bigl|\int_0^1 k_\si(1,y) \sin\la y\,dy\Bigr|
        \to 0,\\
    |\Phi'(\la) - \cos\la| &= \Bigl|\int_0^1 yk_\si(1,y) \cos\la y\,dy\Bigr|
        \to 0
\end{align*}
as $\la\to\infty$ inside the strip
    $\{z\in\bC\mid|\Im z|\le1\}$.
Therefore there exist $\eps_0>0$ and a sequence $(C_n)$ of circles
$C_n:=\{z\in\bC \mid |z-\la_n|= r_n\}$, $r_n\in(0,1/2)$, such that the
discs $D_n:=\{z\in\bC \mid |z-\la_n| \le r_n\}$ are pairwise disjoint and
for all $n\in\bN$ we have
\begin{equation}\label{eq:6.1}
    \min_{\la\in C_n}|\Phi(\la)| \ge \eps_0, \qquad
    \min_{\la\in D_n}|\Phi'(\la)| \ge \eps_0.
\end{equation}

It is proved in~\cite{HM} that the mapping
 \(
    \Si^+ \ni \si \mapsto k_\si(1,\cdot) \in \cH
 \)
is continuous, whence for any $\eps\in(0,\eps_0)$ there exists $\de>0$ such
that for every $\ti\si\in\Si^+$ from a $\de$-neighbourhood of $\si$ we have
$\|k_{\ti\si}(1,\cdot)-k_{\si}(1,\cdot)\|<\eps/4$. We fix such $\eps$ and a
function $\ti\si$, denote by $\ti\la_n^2$ the corresponding eigenvalues and
put~$\ti\Phi(\la):=\Phi(\la,\ti\si)$ and
 $g:=k_{\ti\si}(1,\cdot)-k_\si(1,\cdot)$. Simple computation shows that
\begin{align}\label{eq:6.2}
    \min_{\la\in C_n}|\ti\Phi(\la)-\Phi(\la)| &\le
        2\int_0^1|g(y)|\,dy\le\eps_0/2,\\
    \min_{\la\in D_n}|\ti\Phi'(\la)-\Phi'(\la)| &\le
        2\int_0^1|y g(y)|\,dy\le\eps_0/2.\label{eq:6.3}
\end{align}
Estimates~\eqref{eq:6.1} and \eqref{eq:6.2} and Rouch\'e's theorem imply
that $\ti\la_n\in D_n$.

Since the functions $\ti\Phi$ and $\Phi$ assume real values for real $\la$,
for every $n\in\bN$ there is $\widehat\la_n$ between $\la_n$ and
$\ti\la_n$, for which
\begin{equation}
    \ti\Phi(\la_n) = \ti\Phi(\la_n) - \ti \Phi(\ti\la_n)
        = (\la_n-\ti\la_n) \Phi'(\widehat\la_n).
\end{equation}
We observe that $\widehat\la_n\in D_n$ and that in view of~\eqref{eq:6.1}
and \eqref{eq:6.3},
\[
    |\ti\Phi'(\widehat\la_n)|\ge |\Phi(\widehat\la_n)|-\eps_0/2\ge\eps_0/2.
\]
On the other hand,
\[
    \ti\Phi(\la_n) = \ti\Phi(\la_n) - \Phi(\la_n)
        = s_n(g):= \int_0^1 g(y)\sin\la_ny\,dy,
\]
so that we get
\[
    |\la_n-\ti\la_n| \le \frac{2 s_n(g)}{\eps_0}
\]
for all $n\in\bN$. Recall (see Appendix~\ref{sec:bases}) that the system of
functions $\{\sin\la_nx\}$ forms a Riesz basis of $\cH$, so that
\[
    \sum_{n\in\bN}|s_n(g)|^2 \le C\|g\|^2,
\]
$C>0$ being the corresponding Riesz constant. Combining the above
estimates, we conclude that
\[
    \sum_{n\in\bN} |\la_n-\ti\la_n|^2 \le \frac{4C}{\eps^2_0}\|g\|^2
        \le \frac{4C\eps^2}{\eps^2_0},
\]
and the desired continuity follows.
\end{proof}

\begin{lemma}
The sequence $(\al_k)$ of norming constants depends continuously on
$\hat\si\in\Si^+/\bR$.
\end{lemma}

\begin{proof}
Suppose that $\ti\si\in\cH$ is in the $\eps$-neighbourhood $\cO_\eps$
of~$\si$ and $\ti\al_k$ are the corresponding norming constants for the
operator $T_{\ti\si}$. As in the proof of Lemma~\ref{lem:alk} we see that
\[
    \al_k-\ti\al_k = \|u_k-\ti u_k\| \|u_k+\ti u_k\|,
\]
where $u_k = (I+K_\si)v_k$ and $\ti u_k = (I+ K_{\ti\si})\ti v_k$, $K_\si$
and $K_{\ti\si}$ are transformation operators, and $v_k=\sqrt{2}\sin\la_k
x$, $\ti v_k=\sqrt{2}\sin\ti\la_k x$. Since the norms $\|\ti u_k\|$ are
bounded uniformly in $\ti\si\in\cO_\eps$, it remains to show that the
sequence $\|u_k-\ti u_k\|$ is in $\ell_2$ and
 \(
    \sum_{k\in\bN}\|u_k-\ti u_k\|^2 \to 0
 \)
as $\eps \to 0$.

We have
\[
    u_k - \ti u_k = v_k - \ti v_k + (K_\si - K_{\ti\si})v_k +
        K_{\ti\si}(v_k-\ti v_k).
\]
The system $(v_k)$ with $\la_k=\pi k +\mu_k$ and $(\mu_k)\in\ell_2$ depends
continuously on $(\mu_k)$ in the sense that
\[
    \sum_{k\in\bN}\|v_k - \ti v_k\|^2 \le 2
    \sum_{k\in\bN}\|\mu_k-\ti\mu_k\|^2.
\]
Since the norms $\|K_{\ti\si}\|$ are uniformly bounded in $\ti\si\in
O_\eps$, also
\[
    \sum_{k\in\bN}\|K_{\ti\si}(v_k - \ti v_k)\|^2 \le C_1
    \sum_{k\in\bN}\|\mu_k-\ti\mu_k\|^2.
\]
Finally, we observe that the system $(v_k)$ is a Riesz basis of~$\cH$, so
that
\[
    \sum_{k\in\bN}\|(K_\si - K_{\ti\si})v_k\|^2 \le C_2
        \|K_\si - K_{\ti\si}\|^2_{\fS_2},
\]
see Appendix~\ref{sec:HS}. Since $K_\si$ depends continuously in $\fS_2$ on
$\si\in\cH$ by the results of~\cite{HM}, the required continuity follows.
\end{proof}

Combining the above lemmata, we arrive at the following result.

\begin{theorem}\label{thm:homeom}
The mapping
\[
    \Si^+ / \bR \ni \hat\si \mapsto \{(\la_k^2),(\al_k)\} \in \SD
\]
is a homeomorphism.
\end{theorem}

Fixing the spectrum $(\la_k^2)$, we can say even more about the
corresponding isospectral set (cf.~\cite[Ch.~4]{PT} for the regular case
$q\in\cH$).

\begin{theorem}\label{thm:isosp}
Suppose that the sequence $(\la^2_k)$ satisfies assumption (A1). Then the
set of all isospectral potentials in $W^{-1}_2(0,1)$ with the Dirichlet
spectrum~$(\la_k^2)$ is analytically diffeomorphic to an open subset of the
Hilbert space~$\ell_2$. The diffeomorphism is performed through the
sequence $(\be_k)$, where $\be_k= \al_k-1$.
\end{theorem}

\begin{proof}
It suffices to note that the correspondence $(\ti\be_k)\mapsto \phi$,
\[
    \phi(s) = \sum_{k\in\bN} \bigl( \cos\pi ks - \cos\la_ks\bigr)
          - \sum_{k\in\bN} \ti\be_k \cos\la_ks,
\]
where $\ti\be_k = \be_k/(1+\be_k)$, is bounded, affine (thus analytic) in
$(\ti\be_k)\in\ell_2$. Since $\al_k= 1+\be_k$ are uniformly bounded away
from zero, the mapping $(\be_k)\mapsto \phi$ is analytic in
$(\be_k)\in\ell_2$, and the result follows from Lemmata~\ref{lem:Kcont} and
\ref{lem:sicont}.
\end{proof}

Analogous stability and isospectrality results hold true for other boundary
conditions considered in the next section.

%
\section{The case of other boundary conditions}%
                \label{sec:third}
%

The solution of the inverse spectral problem presented in
Sections~\ref{sec:cnctn}--\ref{sec:inv} can easily be adapted to the case
of other types of boundary conditions, e.~g., the boundary conditions of
the third type at both endpoints, Dirichlet-Neumann, or Neumann-Dirichlet
ones. Below, we shall briefly discuss the modifications to be made for
these cases.

For $\si \in \cH$ and $H,h\in\bC$, we consider a Sturm-Liouville operator
$T_{\si,H,h}$ given by the differential expression $l_\si$
of~\eqref{eq:Sact} and the boundary conditions
\[
        u^{[1]}(0) - H u(0)=0, \qquad u^{[1]}(1) + h u(1)=0.
\]
(We recall that $u^{[1]}(x) = u'(x) - \si (x) u(x)$ is the quasi-derivative
of $u$.) More precisely, $T_{\si,H,h}$ is given by
\[
    T_{\si,H,h}u = l_\si (u) := - (u^{[1]})' - \si u'
\]
on the domain
\begin{multline*}
    \fD(T_{\si,H,h}) = \{ u \in W^1_1(0,1) \mid u^{[1]} \in W^1_1(0,1),
        \ l_{\si}(u)\in\cH, \\ u^{[1]}(0)-H u(0) = u^{[1]}(1)+h u(1) =0\}.
\end{multline*}
Observe that $T_{\si,H,h}=T_{\si+H,0,h-H}$, so that without loss of
generality we may (and will) assume that $H=0$.

It is known~\cite{SS} that for all real-valued $\si\in\cH$ and $h\in\bR$
the operator $T_{\si,0,h}$ is selfadjoint, bounded below, and has discrete
simple spectrum $(\la_k^2)$, $k\in \bN$. As earlier, upon adding a suitable
constant to the potential $q = \si'$, we can make all the eigenvalues
positive. Denote by $u_k$ the eigenfunction of $T_{\si,0,h}$ corresponding
to the eigenvalue $\la_k^2$ and normalized in such a way that
 $u_k(0)=\sqrt2$, and put $\al_k := \|u_k\|^2$.

Our aim is to solve the inverse spectral problem for $T_{\si,0,h}$, i.e.,
firstly, to describe the set $\SD$ of all {\em spectral data}
$\{(\la_k^2),(\al_k)\}$ that can be obtained by varying real-valued
$\si\in\cH$ and $h\in\bR$ and, secondly, for given spectral data
$\{(\la_k^2),(\al_k)\} \in \SD$, to find the corresponding operator
$T_{\si,0,h}$ (i.e., to find the corresponding primitive $\si\in\cH$ of the
potential $q$ and the number $h\in\bR$).

\begin{lemma}\label{lem:lak'}
Suppose that $\si\in\cH$ is real-valued, $h\in\bR$, and that
$\la_1^2<\la_2^2<\dots$ are eigenvalues of the operator $T_{\si,0,h}$ with
$\la_1^2>0$. Then $\la_k = \pi (k-1) + \mu_k$, where $(\mu_k)\in \ell_2$.
\end{lemma}

\begin{proof}
Denote by $\ti T_{\si,0}$ the extension of the operator $T_{\si,0,h}$
obtained by omitting the boundary condition at the point $x=1$.
By~\cite{HM}, the operators $\ti T_{\si,0}$ and $\ti T_{0,0}$ are similar,
and the similarity is performed by the transformation operator
$I+K_{\si,0}$; here $K_{\si,0}$ is an integral operator of Volterra type,
$(K_{\si,0} u)(x) = \int_0^x k_{\si,0}(x,y) u(y)\,dy$, and the kernel
$k_{\si,0}$ has the property that $k_{\si,0}(x,\cdot)$ is an $\cH$-function
for every $x\in[0,1]$.

Put $u(\,\cdot\,,\la) = (I+K_{\si,0}) v(\,\cdot\,,\la)$, where
 $v(x,\la)=\sqrt{2}\cos\la x$.
Then $u_k = u(\cdot,\la_k)$ and the numbers $\pm\la_k$ are solutions of the
equation $u^{[1]}(1,\la)+h u(1,\la)=0$. Using the properties of the
transformation operator $K_{\si,0}$ (see~\cite[Remark~4.3]{HM}), we find
that
\[
    u^{[1]}(1,\la) = -\sqrt2\la\sin\la
        -\sqrt2\la \int_0^1 g_1(x)\sin\la x\, dx
        + \int_0^1 g_2(x)\sqrt2\cos\la x\, dx + C
\]
for some functions $g_1,g_2$ from $\cH$ and a real constant $C$. It follows
that $\la_k$ are zeros of the analytic function
\[
    \Phi_1(\la) := -\la\sin\la + h \cos\la
        -\la \int_0^1 g_1(x)\sin\la x\, dx
        + \int_0^1 g_3(x)\cos\la x\, dx + C/\sqrt2,
\]
where $g_3:=g_2 + h k_{\si,0}(1,\cdot) \in\cH$. Now the standard analysis
(cf. the proof of Lemma~\ref{lem:lak}) yields the asymptotics required.
\end{proof}

We next establish the asymptotics of the norming constants $\al_k$.

\begin{lemma}\label{lem:alk'}
Suppose that $\si\in\cH$ is real-valued, $h\in\bR$, and that $u_k$ are
eigenfunctions of the operator $T_{\si,0,h}$ normalized as above. Then
$\al_k = 1 + \be_k$, where the sequence $(\be_k)$ belongs to~$\ell_2$.
\end{lemma}

Proof of this statement is completely analogous to the proof of
Lemma~\ref{lem:alk}. The minor changes to be made concern notations,
namely, we should put $v_k(s) = v(s,\la_k)= \sqrt2\cos \la_k s$ and
$v_{k,0}(s) = \sqrt2\cos[\pi (k-1)s]$.

Suppose that the sequences $(\la_k^2)$ and $(\al_k)$ are the spectral data
for an operator $T_{\si,0,h}$ with $\si\in\cH$ and $h\in\bR$ and put
\[
    F := \slim_{N\to\infty}\sum_{k=1}^N
        \Bigl(\frac1\al_k (\cdot, v_k)v_k
        -  (\cdot,v_{k,0}) v_{k,0}\Bigr).
\]
It is easily seen that $I + F$ is uniformly positive. As in
Lemma~\ref{lem:F}, it can be proved that $F$ has the following properties.

\begin{lemma}\label{lem:7.3}
The operator $F$ is a Hilbert-Schmidt integral operator with kernel
\begin{equation}\label{eq:F'}
    f(x,y) = \phi(x+y) + \phi(x-y)
\end{equation}
 where
\begin{equation}\label{eq:phi'}
    \phi(s) = \sum_{k=1}^\infty
        \Bigl(\frac1\al_k \cos\la_ks -\cos\pi(k-1)s \Bigr)
\end{equation}
is an $L_2(0,2)$-function.
\end{lemma}

Next we show that $F$ is related to the transformation operator $I +
K_{\si,0}$ through the GLM equation~\eqref{eq:KGLM}, which can be used to
determine uniquely $K_{\si,0}$ through the formula
    $K_{\si,0} = - (\cI + \cP^+_F)^{-1} \cP^+ F \in \fS_2^+$,
see Section~\ref{sec:GLM}.

Suppose now that the sequences $(\la_k^2)$ and $(\al_k)$ consist of
positive numbers and are as claimed in Lemmata~\ref{lem:lak'} and
\ref{lem:alk'}. We solve the GLM equation~\eqref{eq:KGLM} for $K$, and then
in the same manner as in Theorem~\ref{thm:findK} we show that there exists
$\si\in\cH$ such that $I+K$ is the transformation operator for the
Sturm-Liouville operators $\cT_{\si,0}$ and $\cT_{0,0}$. The only thing
that remains to be proved is that our initial sequences $(\la_k^2)$ and
$(\al_k)$ are the spectral data for the Sturm-Liouville operator
$T_{\si,0,h}$ with some $h\in\bR$.

We use identity~\eqref{eq:lagr} and replicate the arguments of the proof of
Theorem~\ref{thm:spdt} to show that none of $u_j(1)$ vanishes as otherwise
all $u_j(1)$ would vanish and $\la_k$ would have a different asymptotics
(that of the Neumann-Dirichlet boundary conditions case, see below).
Therefore $u_j(1)\ne0$ for all $j\in\bN$, and relation~\eqref{eq:lagr}
implies that there exists a constant $h \in \bR$ such that
\begin{equation}\label{eq:h}
    u_j^{[1]}(1)/u_j(1) = h
\end{equation}
for all $j \in \bN$. Thus $(\la^2_k)$ are eigenvalues of the
Sturm-Liouville operator \( -\dd + \si' \) subject to the boundary
conditions $y^{[1]}(0)=0$ and $y^{[1]}(1) = h y(1)$, and the proof is
complete.

We summarize the above considerations in the following theorem.


\begin{theorem}\label{thm:DN}
For two sequences $(\la_k^2)$ and $(\al_k)$ to be the spectral data of a
(positive) Sturm-Liouville operator $T_{\si,0,h}$ with real-valued
potential $q = \si'$ from $W^{-1}_2(0,1)$ and $h\in\bR$, it is necessary
and sufficient that $\al_k$ are as in (A2) and $\la_k$ satisfy the
following assumption:
\begin{itemize}
\item[(A1$'$)] $\la_k$ are all positive, strictly increase with $k$, and
    obey the asymptotic relation $\la_k = \pi (k-1) + \mu_k$
    with some $\ell_2$-sequence $(\mu_k)_{k=1}^\infty$.
\end{itemize}

The corresponding operator $T_{\si,0,h}$ is uniquely recovered from the
spectral data through formula~\eqref{eq:sigma}, in which the functions
$\phi$ and $f$ are given by~\eqref{eq:phi'} and \eqref{eq:F'} respectively,
and $k$ is the kernel of the solution $K$ of the GLM
equation~\eqref{eq:KGLM}. The number $h$ in the boundary conditions is
given by~\eqref{eq:h}.
\end{theorem}

In a similar manner the case where one of the boundary conditions is a
Dirichlet one and the other one is of the third type can be treated. It
suffices to consider only the cases where $H=\infty$, $h=0$, or $H=0$,
$h=\infty$ as other situations reduce to one of these in view of the
relation $T_{\si,H,h}=T_{\si+h',H-h',h+h'}$.

The operators $T_{\si,\infty,0}$ and $T_{\si,0,\infty}$ can be uniquely
recovered from the spectral data, the sequences of eigenvalues and norming
constants. The eigenvalues $\la_k^2$ of the operators $T_{\si,\infty,0}$
and $T_{\si,0,\infty}$ obey the asymptotics $\la_k = \pi (k-1/2) + \mu_k$
for $\ell_2$-sequences $(\mu_k)$, while the norming constants $\al_k$
(defined as in Section~\ref{sec:direct} or Section~\ref{sec:third}
according as $H=\infty$ or $H=0$) satisfy (A2). Thus the sets of spectral
data for the families of Sturm-Liouville operators $T_{\si,\infty,0}$ and
$T_{\si,0,\infty}$, where $\si$ runs through~$\cH$, admit explicit
descriptions, while the reconstruction algorithm remains the same. The
details can be recovered by analogy with the analysis of
Sections~\ref{sec:direct}--\ref{sec:third} (cf.~\cite[Sect.~II.10]{Le} for
the regular case).

%
                        \appendix
%
%
\section{Riesz bases}\label{sec:bases}
%

In this appendix we gather some well known facts about Riesz bases of sines
and cosines (see, e.g., \cite{GK1,HV} and references therein for a detailed
exposition of this topic).

Recall that a sequence $(e_n)_{n=1}^\infty$ in a Hilbert space $\cH$ is a
\emph{Riesz basis} if and only if any element $e\in\cH$ has a unique
expansion
        \( e = \sum_{n=1}^\infty c_n e_n\)
with $(c_n)\in \ell_2$. If $(e_n)$ is a Riesz basis, then in the above
expansion the \emph{Fourier coefficients} $c_n$ are given by $c_n =
(e,e'_n)$, where $(e'_n)_{n=1}^\infty$ is a system biorthogonal to $(e_n)$,
i.e., which satisfies the equalities $(e_k,e'_n) = \de_{kn}$ for all
$k,n\in\bN$. Moreover, the biorthogonal system $(e'_n)$ is a Riesz basis of
$\cH$ as long as $(e_n)$ is, in which case for any $e\in\cH$ also the
expansion $e=\sum(e,e_n)e'_n$ takes place. In particular, if $(e_n)$ is a
Riesz basis, then for any $e\in\cH$ the sequence $(c'_n)$ with
$c'_n:=(e,e_n)$ belongs to $\ell_2$.

Also, any Riesz basis $(e_n)$ is equivalent to an orthonormal one, i.e.,
there exists a homeomorphism $U$ such that $(Ue_n)$ is an orthonormal basis
of~$\cH$. As a result, there exists a constant $C>0$ (the {\em Riesz
constant}) such that, for any sequence $(c_k)\in\ell_2$,
\[
    C^{-1} \sum |c_k|^2 \le
        \|\sum c_k e_k\|^2 \le
            C \sum |c_k|^2;
\]
in particular, the series $\sum c_k e_k$ converges in $\cH$ for any
$\ell_2$-sequence~$(c_k)$.

\begin{proposition}[\cite{HV}]\label{pro:RB}
Suppose that the real numbers $\mu_k$, $k\in\bN$, tend to zero and that the
sequence $(\pi k+\mu_k)$ strictly increases. Then each of the following
systems forms a Riesz basis of $L_2(0,1)$:
\begin{itemize}
\item [(a)] $\{\sin(\pi kx + \mu_kx)\}_{k=1}^\infty$;
\item [(b)] $\{\sin(\pi [k-1/2]x + \mu_kx)\}_{k=1}^\infty$;
\item [(c)] $\{\cos(\pi kx + \mu_kx)\}_{k=0}^\infty$;
\item [(d)] $\{\cos(\pi [k+1/2]x + \mu_kx)\}_{k=0}^\infty$.
\end{itemize}
\end{proposition}


\begin{corollary}\label{cor:RB}
Under the assumptions of Proposition~\ref{pro:RB} the systems~(a)--(d)
constitute Riesz bases of their closed linear spans in $L_2(0,2)$.
\end{corollary}


%
\section{The ideal of Hilbert-Schmidt operators}\label{sec:HS}
%

In this appendix we recall some necessary facts about the Schatten-von
Neumann $\fS_p$ ideals (see details in~\cite{GK1}).

Suppose that $T$ is a compact operator in a Hilbert space $\cH$;
then $|T| = (T^*T)^{1/2}$ is a nonnegative selfadjoint compact operator.
Denote by $\la_1(|T|) \ge \la_2(|T|) \ge \dots$ the eigenvalues of $|T|$ in
non-decreasing order and repeated according to their multiplicity.
Then $\la_k(|T|)$ is called the {\em $k$-th $s$-number} of $T$ and
is denoted by $s_k(T)$.

The ideal $\fS_p$, $p\in[1,\infty)$, consists of all compact operators for
which the expression
\[
    \|T\|_{\fS_p} := \Bigl( \sum s_k^p(T)\Bigr)^{1/p}
\]
is finite. Being endowed with the norm $\|\,\cdot\,\|_{\fS_p}$, the ideal
$\fS_p$ becomes a Banach space.

In particular, $\fS_1$ is the ideal of trace class operators. For any
$T\in\fS_1$ its {\em matrix trace} $\tr T := \sum_{k\in\bN} (T e_k,e_k)$
with respect to any orthonormal basis $(e_k)$ coincides with the {\em
spectral trace}, i.e., the sum of all eigenvalues of $T$ repeated according
to their algebraic multiplicity, and is finite.

For $p=2$ the ideal $\fS_2$ consists of all Hilbert-Schmidt operators. For
$T\in \fS_2$ and any two orthonormal bases $(e_j)$ and $(e'_k)$ of $\cH$,
we have
\[
    \|T\|_{\fS_2}^2 = \sum_{j,k\in\bN} |(T e_j,e'_k)|^2
            = \sum_{j\in\bN} \|T e_j\|^2 < \infty.
\]
If $(e_j)$ is a Riesz basis of $\cH$, then $(Ue_j)$ is an orthonormal basis
for some homeomorphism~$U$, and for any $T\in\fS_2$ the estimate
\[
     \sum_{j\in\bN} \|T e_j\|^2 = \|TU^{-1}\|_{\fS_2}^2
        \le\|T\|_{\fS_2}^2 \|U^{-1}\|^2
\]
holds. Also, $ST$ is a trace class operator whenever $S$ and $T$ are
Hilbert-Schmidt ones. Moreover, $\fS_2$ is a Hilbert space with the scalar
product given by
\[
    \lan S,T\ran_2 := \tr S T^* = \tr T^* S.
\]
If $\cH = L_2(0,1)$, any Hilbert-Schmidt operator $T$ is an integral one
with kernel $t$ given by
\[
    t(x,y):= \sum_{j,k\in\bN} (Te_j,e_k)e_k(x)\ov{e_j(y)};
\]
here $(e_k)$ is any orthonormal basis of $\cH$ and the series converges in
$L_2\bigl((0,1)\times(0,1)\bigr)$. Moreover,
\[
    \int_0^1 \int_0^1 |t(x,y)|^2 \,dx \,dy = \|T\|_{\fS_2}^2.
\]


\end{document}